# Stabilizing the Consistent Quasidiffusion Method with Linear Prolongation

Dean Wang[1,*]

[1]The Ohio State University, Columbus, Ohio

## ABSTRACT

The quasidiffusion (QD) method, also known as the Variable Eddington Factor (VEF) method in the astrophysical community, is an established iterative method for accelerating source iterations in $S_N$ calculations. A great advantage of the QD method is that the diffusion equation that accelerates the $S_N$ source iterations can be discretized in any valid discretization without concern for consistency with the transport discretization. QD has comparable effectiveness with diffusion synthetic acceleration (DSA), but the converged scalar flux of the diffusion equation will differ from the transport solution by the spatial truncation errors. Larsen et al. introduced a new consistent QD method (CQD), which includes a straightforwardly defined transport consistency factor closely related to the well-known coarse mesh finite difference (CMFD) and DSA methods. The CQD method preserves the discretized scalar flux solution of the $S_N$ equations, and it is stable for problems with optically thin spatial cells, but just like nonlinear diffusion acceleration (NDA), it degrades in performance and eventually becomes unstable when the spatial cells become greater than about one mean free path thick. In this paper, we performed a formal Fourier analysis of the CQD method to show that its theoretical spectral radius is essentially the same as that of the NDA method. To improve the stability of CQD, we introduce the lpCQD method, which adopts the idea of the linear prolongation CMFD (lpCMFD) method.

## 1.     INTRODUCTION

The QD method (or VEF) was one of the first nonlinear methods for accelerating source iterations (SIs) in $S_N$ calculations [1, 2]. It is comparable in effectiveness to both linear and nonlinear forms of DSA [3], but it offers much more flexibility than DSA. Stability can only be guaranteed with DSA if the diffusion equation is differenced in a manner consistent with that of the $S_N$ equations. Larsen et al. modified the QD method to make it a true acceleration scheme by including a straightforwardly defined transport consistency factor in the discretized low-order QD equations, which is closely related to the well-known CMFD (or NDA) and DSA methods [4]. Compared to the conventional QD equations, the new QD method, called CQD, is simpler since it does not require the calculation and storage of cell-edge scalar flux estimates. CQD preserves the discretized scalar flux solution of the $S_N$ equations, and it is stable for problems with optically thin spatial cells, but just like CMFD, it degrades in performance and eventually becomes unstable when the spatial cells become greater than about one mean free path thick. It should be noted that there exist early works to modify QD to make it consistent with the transport solution [5, 6].

---

*wang.12239@osu.edu

In this paper, we perform a formal Fourier analysis of the CQD method and show that the theoretical spectral radius of CQD is the same as that of NDA. To improve the stability of CQD, we apply the "linear prolongation" idea of the lpCMFD scheme [7] to the CQD method. This new CQD method is called "lpCQD".

The remainder of this paper is organized as follows. Section 2 presents the Fourier analysis of the CQD method to explain why the CQD has very similar stability properties to NDA. Section 3 introduces the linear prolongation approach to stabilize CQD. Numerical results are given in Section 4. Section 5 contains a summary and discussion.

## 2. FOURIER ANALYSIS OF CQD

Our model problem is a one-dimensional (1-D) slab geometry, homogeneous monoenergetic problem with isotropic scattering and periodic boundary conditions. This model problem is solved with the neutron transport equation coupled with CQD acceleration.

### 2.1. Algorithm of CQD

Using standard notation, the coupled iteration procedure is described as follows.

1) *Perform a transport sweep:*

$$\mu \frac{\partial \psi^{l+\frac{1}{2}}(x,\mu)}{\partial x} + \Sigma_t(x)\psi^{l+\frac{1}{2}}(x,\mu) = \frac{1}{2}\Sigma_s(x)\phi^l(x) + \frac{1}{2}Q(x), \tag{1}$$

where

$$\phi(x) = \int_{-1}^{1} \psi(x,\mu)d\mu. \tag{2}$$

2) *Compute the Eddington factor and the nonlinear current correction factor for the low-order diffusion equation:*

We define the diffusion coefficient and the nonlinear current correction factor as

$$E^{l+\frac{1}{2}} = \frac{\phi_{(2)}^{l+1/2}}{\phi^{l+\frac{1}{2}}}, \tag{3a}$$

$$\widetilde{D}^{l+\frac{1}{2}} = \frac{J^{l+\frac{1}{2}} + \frac{1}{\Sigma_t}\frac{d}{dx}\phi_{(2)}^{l+\frac{1}{2}}}{\phi^{l+\frac{1}{2}}}, \tag{3b}$$

where

$$\phi_{(2)}(x) = \int_{-1}^{1} \psi(x,\mu)\mu^2 d\mu. \tag{3c}$$

3) *Solve the low-order diffusion problem:*

$$\frac{d}{dx}\left(-\frac{1}{\Sigma_t}\frac{d}{dx}E^{l+\frac{1}{2}}\phi^{l+1} + \widetilde{D}^{l+\frac{1}{2}}\phi^{l+1}\right) + \Sigma_a\phi^{l+1} = Q(x). \tag{4}$$

4) *Update the scalar flux:*

$$\phi^{l+1} = \phi^{l+\frac{1}{2}} + \delta\phi^{l+1}, \tag{5a}$$

where $\delta\phi^{l+1}$ is a genera representation of different flux updating schemes, which will be given later. For the standard flux updating method, e.g., NDA (or CMFD):

$$\delta\phi^{l+1} = \phi^{l+1} - \phi^{l+\frac{1}{2}}. \tag{5b}$$

## 2.2. Linearization of CQD

Fourier analysis is typically performed based on a 1-D homogeneous problem with a uniform mesh. The monoenergetic $S_N$ transport equation is firstly linearized near the exact solution. Then the error terms in the linearized equations are expressed as Fourier modes. Fourier analysis is carried out on a single mesh cell with periodic boundary conditions. The mesh spacing is $h$. The spectral radius can be obtained from the error iteration matrix.

We follow the standard approach to perform Fourier analysis. First, we define the iterative unknowns of the model problem:

$$\psi(x,\mu) = \frac{\Phi_0}{2} + \varepsilon\psi_1(x,\mu), \tag{6a}$$

$$\phi(x) = \Phi_0 + \varepsilon\phi_1(x), \tag{6b}$$

$$J(x) = \int_{-1}^{1} \psi(x,\mu)\mu d\mu = \varepsilon \int_{-1}^{1} \psi_1(x,\mu)\mu d\mu = \varepsilon J_1(x), \tag{6c}$$

$$\phi_{(2)} = \int_{-1}^{1} \psi(x,\mu)\mu^2 d\mu = \int_{-1}^{1}\left(\frac{\Phi_0}{2} + \varepsilon\psi_1(x,\mu)\right)\mu^2 d\mu = \frac{\Phi_0}{3} + \varepsilon\phi_{(2),1}, \tag{6d}$$

where

$$\phi_{(2),1} = \int_{-1}^{1} \psi_1(x,\mu)\mu^2 d\mu. \tag{6e}$$

Now we linearize the neutron transport equation. Substituting Eqs. (6a) and (6b) into (1), after some algebra we obtain by neglecting the $O(\varepsilon^2)$ terms:

$$\mu\frac{\partial}{\partial x}\psi_1^{l+\frac{1}{2}} + \Sigma_t\psi_1^{l+\frac{1}{2}} = \frac{c}{2}\Sigma_t\phi_1^l. \tag{7}$$

Next, we linearize the diffusion equation, (4). First, we linearize the term $E^{l+\frac{1}{2}}\phi^{l+1}$ as follows.

$$E^{l+\frac{1}{2}}\phi_{(0)}^{l+1} = \frac{\phi_{(2)}^{l+\frac{1}{2}}}{\phi^{l+\frac{1}{2}}}\phi_{(0)}^{l+1} = \left(\frac{\frac{\Phi_0}{3} + \varepsilon\phi_{(2),1}^{l+\frac{1}{2}}}{\Phi_0 + \varepsilon\phi_1^{l+\frac{1}{2}}}\right)(\Phi_0 + \varepsilon\phi_1^{l+1})$$

$$= \left(\frac{\frac{\Phi_0 + \varepsilon\phi_1^{l+\frac{1}{2}}}{3} - \frac{\varepsilon\phi_1^{l+\frac{1}{2}}}{3} + \varepsilon\phi_{(2),1}^{l+\frac{1}{2}}}{\Phi_0 + \varepsilon\phi_1^{l+\frac{1}{2}}}\right)(\Phi_0 + \varepsilon\phi_1^{l+1})$$

$$= \left(\frac{1}{3} + \frac{-\frac{\varepsilon\phi_1^{l+\frac{1}{2}}}{3} + \varepsilon\phi_{(2),1}^{l+\frac{1}{2}}}{\Phi_0 + \varepsilon\phi_1^{l+\frac{1}{2}}}\right)(\Phi_0 + \varepsilon\phi_1^{l+1})$$

$$\approx \frac{1}{3}(\Phi_0 + \varepsilon\phi_1^{l+1}) - \frac{\varepsilon\phi_1^{l+\frac{1}{2}}}{3} + \varepsilon\phi_{(2),1}^{l+\frac{1}{2}}. \tag{8}$$

Then we linearize the term $\widetilde{D}^{l+\frac{1}{2}}\phi^{l+1}$ as:

$$\widetilde{D}^{l+\frac{1}{2}}\phi^{l+1} = \frac{J^{l+\frac{1}{2}} + \frac{1}{\Sigma_t}\frac{d}{dx}\phi_{(2)}^{l+\frac{1}{2}}}{\phi^{l+\frac{1}{2}}} = \frac{\varepsilon J_1^{l+\frac{1}{2}} + \varepsilon\frac{1}{\Sigma_t}\frac{d}{dx}\phi_{(2),1}^{l+\frac{1}{2}}}{\Phi_0 + \varepsilon\phi_1^{l+\frac{1}{2}}}(\Phi_0 + \varepsilon\phi_1^{l+1})$$

$$\approx \varepsilon J_1^{l+\frac{1}{2}} + \varepsilon\frac{1}{\Sigma_t}\frac{d}{dx}\phi_{(2),1}^{l+\frac{1}{2}}, \tag{9}$$

Substituting Eqs. (8) and (9) into (4), we obtain

$$\frac{d}{dx}\left[-\frac{1}{\Sigma_t}\frac{d}{dx}\left(\frac{1}{3}(\Phi_0 + \varepsilon\phi_1^{l+1}) - \frac{\varepsilon\phi_1^{l+\frac{1}{2}}(x)}{3} + \varepsilon\phi_{(2),1}^{l+\frac{1}{2}}\right) + \varepsilon J_1^{l+\frac{1}{2}} + \varepsilon\frac{1}{\Sigma_t}\frac{d}{dx}\phi_{(2),1}^{l+\frac{1}{2}}\right]$$

$$+\Sigma_a(\Phi_0 + \varepsilon\phi_1^{l+1}) = Q. \tag{10}$$

Simplifying the above equation:

$$\frac{1}{3\Sigma_t}\frac{d^2}{dx^2}\phi_1^{l+1}(x) - \Sigma_a\phi_1^{l+1} = \frac{1}{3\Sigma_t}\frac{d^2}{dx^2}\phi_1^{l+\frac{1}{2}}(x) + \frac{d}{dx}J_1^{l+\frac{1}{2}}. \tag{11}$$

Note that the two terms of $\phi_{(2),1}^{l+\frac{1}{2}}$ in Eq. (10) have canceled out each other, and the linearized CQD equation is reduced to that of NDA. Therefore, it should be expected that the CQD method have the same spectral radius as the conventional NDA method with Fourier analysis.

We rewrite Eqs. (7) and (11) by dropping the subscript "1" as:

$$\mu\frac{\partial}{\partial x}\psi^{l+\frac{1}{2}} + \Sigma_t\psi^{l+\frac{1}{2}} = \frac{c}{2}\Sigma_t\phi^l. \tag{12a}$$

$$\frac{1}{3\Sigma_t}\frac{d^2}{dx^2}\phi^{l+1}(x) - \Sigma_a\phi^{l+1} = \frac{1}{3\Sigma_t}\frac{d^2}{dx^2}\phi^{l+\frac{1}{2}}(x) + \frac{d}{dx}J^{l+\frac{1}{2}}. \tag{12b}$$

### 2.3. Fourier Analysis

To perform Fourier analysis, the linearized transport equation, Eq. (12a), is typically discretized using the DD or SC method a mesh grid, as shown in Figure 1.

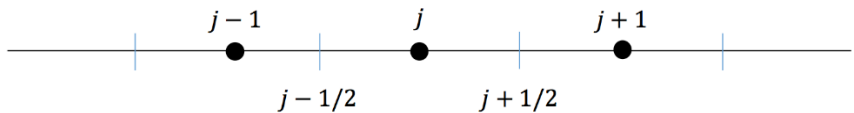

**Figure 1. 1-D mesh.**

$$\frac{\mu_n}{h}\left(\psi_{n,j+\frac{1}{2}}^{l+\frac{1}{2}} - \psi_{n,j-\frac{1}{2}}^{l+\frac{1}{2}}\right) + \Sigma_t \psi_{n,j}^{l+\frac{1}{2}} = \frac{c}{2}\Sigma_t \phi_j^l, \tag{13a}$$

where

$$\psi_{n,j} = \left(\frac{1-\alpha_n}{2}\right)\psi_{n,j-\frac{1}{2}} + \left(\frac{1+\alpha_n}{2}\right)\psi_{n,j+\frac{1}{2}}, \tag{13b}$$

$$\phi_j = \sum_{n=1}^{N} \psi_{n,j} w_n, \tag{13c}$$

$$\alpha_n = \begin{cases} 0, & \text{DD} \\ \dfrac{1+e^{-\frac{\Sigma_t h}{\mu_n}}}{1-e^{-\frac{\Sigma_t h}{\mu_n}}} - \dfrac{2\mu_n}{\Sigma_t h}, & \text{SC} \end{cases}, \tag{13d}$$

Note that the mesh size is assumed to be uniform. SC denotes the step characteristic method.

The diffusion equation, Eq. (12b) is discretized using the central differencing scheme on the same fine-mesh grid as for the $S_N$ transport equation:

$$\frac{1}{3\Sigma_t}\frac{\phi_{j-1}^{l+1} - 2\phi_j^{l+1} + \phi_{j+1}^{l+1}}{h^2} - \Sigma_a \phi_j^{l+1} = \frac{1}{3\Sigma_t}\frac{\phi_{j-1}^{l+\frac{1}{2}} - 2\phi_j^{l+\frac{1}{2}} + \phi_{j+1}^{l+\frac{1}{2}}}{h^2} + \frac{J_{j+\frac{1}{2}}^{l+\frac{1}{2}} - J_{j-\frac{1}{2}}^{l+\frac{1}{2}}}{h}. \tag{14}$$

We introduce the following Fourier ansatz:

$$\phi_j^l = \lambda^l e^{i\Sigma_t \omega x_j}, \tag{15a}$$

$$\psi_{n,j}^{l+\frac{1}{2}} = \lambda^l a_n e^{i\Sigma_t \omega x_j}, \tag{15b}$$

$$\psi_{n,j\pm\frac{1}{2}}^{l+\frac{1}{2}} = \lambda^l b_n e^{i\Sigma_t \omega x_{j\pm\frac{1}{2}}}, \tag{15c}$$

$$\phi_j^{l+\frac{1}{2}} = \lambda^l \beta e^{i\Sigma_t \omega x_j}, \tag{15d}$$

$$\phi_j^{l+1} = \lambda^l \gamma e^{i\Sigma_t \omega x_j}. \tag{15e}$$

We substitute Eqs. (15a) – (15c) into (13a) to obtain

$$2i\frac{\mu_n}{h}\left(\sin\frac{\Sigma_t h\omega}{2}\right) b_n + \Sigma_t a_n = \frac{c}{2}\Sigma_t, \tag{16}$$

and substitute Eqs. (15b) and (15c) into (13b) to have

$$a_n = b_n\left[\left(\frac{1-\alpha_n}{2}\right)e^{-i\Sigma_t \omega h/2} + \left(\frac{1+\alpha_n}{2}\right)e^{i\Sigma_t \omega h/2}\right]$$

$$= b_n\left[\cos\frac{\Sigma_t h\omega}{2} + i\alpha_n \sin\frac{\Sigma_t h\omega}{2}\right]. \tag{17}$$

Then we obtain

$$b_n = \frac{a_n}{\left[\cos\frac{\Sigma_t h\omega}{2} + i\alpha_n \sin\frac{\Sigma_t h\omega}{2}\right]} . \tag{18}$$

Substituting Eq. (18) into (16), we have

$$a_n = \frac{c}{2} \frac{1}{1 + \dfrac{i\mu_n \Lambda}{1 + \dfrac{i\alpha_n \Sigma_t h\Lambda}{2}}} , \tag{19a}$$

where

$$\Lambda \equiv \frac{2}{\Sigma_t h} \tan\frac{\Sigma_t h\omega}{2} . \tag{19b}$$

Substituting (15b) and (15d) into (13c), we have

$$\beta = \sum_{n=1}^{N} a_n w_n . \tag{20}$$

Substituting Eq. (19a) into (20), we obtain the spectral radius of SI as:

$$\varrho_{SI}(\omega) = \beta = \sum_{n=1}^{N} a_n w_n = \frac{c}{2} \sum_{n=1}^{N} \frac{w_n}{1 + \dfrac{i\mu_n \Lambda}{\left[1 + \dfrac{i\alpha_n \Sigma_t h\Lambda}{2}\right]}} . \tag{21}$$

Using (15c), we obtain

$$J_{j\pm\frac{1}{2}}^{l+\frac{1}{2}} = \sum_{n=1}^{N} \mu_n \psi_{n,j\pm\frac{1}{2}}^{l+\frac{1}{2}} w_n = \lambda^l e^{i\Sigma_t \omega x_{j\pm\frac{1}{2}}} \sum_{n=1}^{N} b_n \mu_n w_n . \tag{22}$$

We substitute Eqs. (15d), (15e), and (22) into Eq. (14) to find the spectral radius of the CQD method:

$$\varrho_{CQD}(\omega) = \gamma = \beta - \frac{h\left(2i\sin\frac{\Sigma_t h\omega}{2}\right)\sum_{n=1}^{N} b_n \mu_n w_n + h^2 \Sigma_a \beta}{\frac{2}{3\Sigma_t}[1 - \cos(\Sigma_t \omega h)] + h^2 \Sigma_a} , \tag{23a}$$

or

$$\varrho_{CQD}(\omega) = \varrho_{SI}(\omega) - \frac{h\left(2i\sin\frac{\Sigma_t h\omega}{2}\right)\sum_{n=1}^{N} b_n \mu_n w_n + h^2 \Sigma_a \varrho_{SI}(\omega)}{\frac{2}{3\Sigma_t}[1 - \cos(\Sigma_t \omega h)] + h^2 \Sigma_a} . \tag{23b}$$

For the DD method, the above equation can be simplified as

$$\varrho_{CQD}^{DD}(\omega) = \varrho_{SI}(\omega) - \frac{hc\Lambda\left[\tan\left(\frac{\Sigma_t h\omega}{2}\right)\right]\sum_{n=1}^{N} \dfrac{\mu_n^2 w_n}{1 + \mu_n^2 \Lambda^2} + h^2 \Sigma_a \varrho_{SI}(\omega)}{\frac{2}{3\Sigma_t}[1 - \cos(\Sigma_t \omega h)] + h^2 \Sigma_a} . \tag{24}$$

The spectral radius of the CQD method is given as

$$\rho = \max\left(\text{abs}\left(\varrho_{CQD}(\omega)\right)\right) , \tag{25a}$$

where the Fourier frequency for the periodic boundary conditions is expressed as follows:

$$\omega = \frac{2\pi s}{\Sigma_t L}, \quad s = 1, 2, \ldots, J\left(=\frac{L}{h}\right). \tag{25b}$$

For the reflective boundary conditions, the above discrete Fourier frequency is simply divided by 2.

## 3. lpCQD METHOD

In this section, we present the lpCQD method. Like lpCMFD, it replaces the conventional flat flux ratio-based scaling approach in the standard CMFD method with a linear interpolation of the scalar flux differences at the mesh cell edges between the high-order transport and low-order diffusion solutions.

### 3.1. Linear Prolongation

To fix the idea, we use a 1-D problem as an example. In lpCQD, the scalar flux $\phi_j^{l+1}$ for the interior cells $(2 \leq j \leq J-1)$ is updated using a linear prolongation on three neighboring cells as:

$$\phi_j^{l+1} = \phi_j^{l+\frac{1}{2}} + \frac{1}{2}\left[\frac{h_j}{h_{j-1}+h_j}\delta\phi_{j-1}^{l+1} + \left(\frac{h_{j-1}}{h_{j-1}+h_j} + \frac{h_{j+1}}{h_j+h_{j+1}}\right)\delta\phi_j^{l+1} + \frac{h_j}{h_j+h_{j+1}}\delta\phi_{j+1}^{l+1}\right], \tag{26}$$

with

$$\delta\phi_j^{l+1} = \phi_j^{l+1} - \phi_j^{l+\frac{1}{2}}, \tag{27}$$

where $\phi_j^{l+1}$ is calculated from the low-order QD equation and $\phi_j^{l+\frac{1}{2}}$ is the high-order $S_N$ solution. The superscript "$l$" is the index of SI. Note that for the uniform mesh, i.e., $h_j = h$, Eq. (5) can be simplified as:

$$\phi_j^{l+1} = \phi_j^{l+\frac{1}{2}} + \frac{1}{4}\left(\delta\phi_{j-1}^{l+1} + 2\delta\phi_j^{l+1} + \delta\phi_{j+1}^{l+1}\right). \tag{28}$$

Care should be taken for the boundary cells, $j = 1$ and $J$, because $\delta\phi_0^{l+1}$ and $\delta\phi_{J+1}^{l+1}$ are not available. For the reflective BC, since we can simply set $\delta\phi_0^{l+1} = \delta\phi_1^{l+1}$ or $\delta\phi_{J+1}^{l+1} = \delta\phi_J^{l+1}$. However, for the vacuum BC, we suggest the scalar flux update for the boundary cell, e.g., $j = 1$, as follows:

$$\phi_1^{l+1} = \phi_1^{l+\frac{1}{2}} + \frac{1}{2}\left[\left(\alpha + \frac{h_1}{h_1+h_2}\right)\delta\phi_1^{l+1} + \frac{h_1}{h_1+h_2}\delta\phi_2^{l+1}\right], \tag{29}$$

where $0 \leq \alpha \leq 1$. It is noted that the overall convergence is almost insensitive to the choice of $\alpha$.

### 3.2. Fourier Analysis of lpCQD

As in Sect. II, we first linearize the lpCQD flux updating equation, Eq. (28) to obtain equation:

$$\phi_j^{l+1} = \phi_j^{l+\frac{1}{2}} + \frac{1}{4}\left[\left(\phi_{j-1}^{l+1} - \phi_{j-1}^{l+\frac{1}{2}}\right) + 2\left(\phi_j^{l+1} - \phi_j^{l+\frac{1}{2}}\right) + \left(\phi_{j+1}^{l+1} - \phi_{j+1}^{l+\frac{1}{2}}\right)\right]. \tag{30}$$

Note that the above $\phi$'s are small perturbations of the scalar fluxes.

Substituting Eqs. (15d) and (15e) into (30), we obtain

$$\phi_j^{l+1} = \left[\beta + \frac{1}{2}(1 + \cos(\Sigma_t h\omega))(\gamma - \beta)\right]\lambda^l e^{i\Sigma_t \omega x_j}$$

$$= \left[\beta + \frac{1}{2}(1 + \cos(\Sigma_t h\omega))(\gamma - \beta)\right]\phi_j^l, \qquad (31)$$

where $\beta$ and $\gamma$ are given by Eqs. (20) and (23a), respectively. Note that we have used Eq. (15a), i.e., $\phi_j^l = \lambda^l e^{i\Sigma_t \omega x_j}$ in the above calculation. Then we obtain the spectral radius of lpCQD as

$$\varrho_{\text{lpCQD}}(\omega) = \varrho_{\text{SI}}(\omega) + \frac{1}{2}(1 + \cos(\Sigma_t h\omega))\left(\varrho_{\text{CQD}}(\omega) - \varrho_{\text{SI}}(\omega)\right), \qquad (32)$$

where $\varrho_{\text{SI}}(\omega)$ and $\varrho_{\text{CQD}}(\omega)$ are given by Eqs. (21) and (23b), respectively.

## 4. NUMERICAL RESULTS

We implement SI, CQD, and lpCQD in a 1-D $S_N$ code for two different problems. Numerical results are presented below.

### 4.1. Consistency

Problem 1 is designed to demonstrate that SI, CQD and lpCQD are encoded correctly, and both CQD and lpCQD converge to the $S_N$ solution. The physical system is a slab $0 < x < 30$ cm with the reflective boundary at left and vacuum boundary at right. The system has three subregions, each 10 cm thick. Going from left to right, subregions 1 and 3 have $\Sigma_t = 1.0$ cm$^{-1}$, $\Sigma_s = 0.9$ cm$^{-1}$, and $Q = 1.0$ cm$^{-3}$ sec$^{-1}$. Subregion 2 have $\Sigma_t = 1.0$ cm$^{-1}$, $\Sigma_s = 0.99$ cm$^{-1}$, and $Q = 1.0$ cm$^{-3}$ sec$^{-1}$. We ran the problem using (i) the diamond difference (DD) method, (ii) the standard $S_{10}$ Gauss-Legendre quadrature set, (iii) a spatial grid with uniform spatial cells of thickness $h = 1.0$ cm, and (iv) a convergence criterion of $10^{-10}$. The resulting scalar fluxes are plotted in Figure 2. The numerical values of the scalar fluxes agreed to within the $10^{-10}$ convergence criterion for all methods.

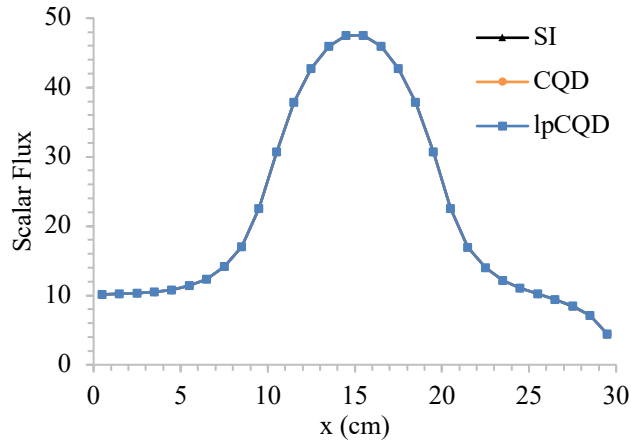

**Figure 2. Scalar fluxes for problem 1.**

### 4.2. Convergence

Problem 2 was designed to assess the effect of linear prolongation on the stability of CQD. Problem 2 is a homogeneous slab of variable thickness $0 \leq x \leq X = 100h$, where $h$ is the (variable) width of a spatial cell. The left boundary is reflective, and the right one is vacuum, and within the system, $\Sigma_t = 1.0$ cm$^{-1}$, $\Sigma_s = 0.4, 0.6, 0.9,$ and $0.99$ cm$^{-1}$, and $Q = 1.0$ cm$^{-3}$ sec$^{-1}$. We used the DD method and the $S_{10}$ quadrature set with a convergence criterion of $10^{-10}$.

Assessing the effect of the spatial cell width $h$ on the rate of convergence of the numerical methods can be done by estimating the numerical spectral radius. To do this, we calculate the $L_2$ norm of the difference between the cell-averaged scalar flux estimates from successive iterates:

$$\|\phi^{l+1} - \phi^l\| = \left(\frac{\sum_{j=1}^{J}|\phi_j^{l+1} - \phi_j^{l+1}|^2 h_j}{\sum_{j=1}^{J} h_j}\right)^{\frac{1}{2}}, \tag{33}$$

and then we estimate the spectral radius $\rho$ by calculating the ratio as follows:

$$\rho \approx \frac{\|\phi^{l+1} - \phi^l\|}{\|\phi^l - \phi^{l-1}\|}. \tag{34}$$

In Figure 3 below, we plot the spectral radius versus the optical cell width $\Sigma_t h$, for SI, CQD, and lpCQD, respectively. It can be seen that the CQD method quickly loses its effectiveness as the optical cell width approaches one and eventually fails, whereas the linear prolongation can not only stabilize the CQD method, but also greatly enhance its acceleration performance for relatively large optical cell width.

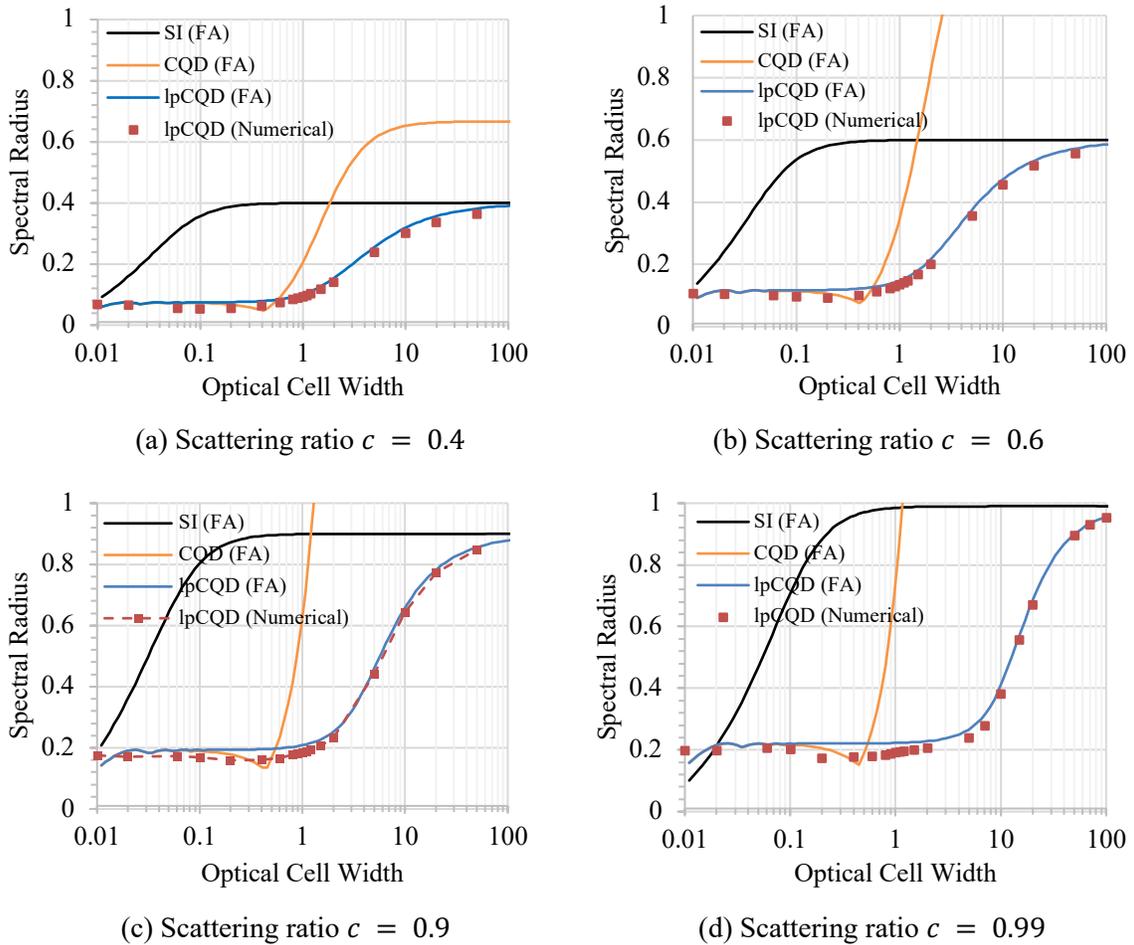

(a) Scattering ratio $c = 0.4$

(b) Scattering ratio $c = 0.6$

(c) Scattering ratio $c = 0.9$

(d) Scattering ratio $c = 0.99$

**Figure 3. Convergence comparison.**

## 5. CONCLUSIONS

We have presented the Fourier analysis to provide a theoretical understanding of the stability of the CQD method. It has shown that the linearized equations of the CQD method are the same as those of the conventional NDA method, and thus they share the similar convergence performance and stability characteristics. To enhance the CQD method, we have proposed the lpCQD method, which employs the linear prolongation approach for flux updating, similar to the lpCMFD method. The Fourier analysis and numerical results demonstrate that the new lpCQD method exhibits better convergence and stability than the CQD method.